
\documentstyle{amsppt}
\baselineskip18pt
\magnification=\magstep1
\pagewidth{30pc}
\pageheight{45pc}
\hyphenation{co-deter-min-ant co-deter-min-ants pa-ra-met-rised
pre-print pro-pa-gat-ing pro-pa-gate
fel-low-ship Cox-et-er dis-trib-ut-ive}
\def\leaderfill{\leaders\hbox to 1em{\hss.\hss}\hfill}
\def\A{{\Cal A}}
\def\D{{\Cal D}}
\def\H{{\Cal H}}

\

\def\afn{{\text {\bf a}}}

\def\idest{i.e.,\ }

\def\be{{\beta}}
\def\g{{\gamma}}

\def\d{{\delta}}

\def\e{{\varepsilon}}

\def\th{{\theta}}

\def\k{{\kappa}}
\def\l{{\lambda}}

\def\s{{\sigma}}
\def\t{{\tau}}

\def\bc{{\bold c}}

\def\BB{{\bold B}}
\def\b0{\text{\bf 0}}

\def\ra{{\ \longrightarrow \ }}

\def\lan{{\langle}}
\def\ran{{\rangle}}

\def\real{{\Bbb R}}
\def\complex{{\Bbb C}}
\def\zed{{\Bbb Z}}

\def\Im{\text{\rm Im}}
\def\End{\text{\rm End}}

\def\boxit#1{\vbox{\hrule\hbox{\vrule \kern3pt
\vbox{\kern3pt\hbox{#1}\kern3pt}\kern3pt\vrule}\hrule}}
\def\rabbit{\vbox{\hbox{\kern0pt
\vbox{\kern0pt{\hbox{---}}\kern3.5pt}}}}

\def\tableau#1{
        \hbox {
                \hskip -10pt plus0pt minus0pt
                \raise\baselineskip\hbox{
                \offinterlineskip
                \hbox{#1}}
                \hskip0.25em
        }
}

\def\tabCol#1{
\hbox{\vtop{\hrule
\halign{\strut\vrule\hskip0.5em##\hskip0.5em\hfill\vrule\cr\lower0pt
\hbox\bgroup$#1$\egroup \cr}
\hrule
} } \hskip -10.5pt plus0pt minus0pt}

\def\CR{
        $\egroup\cr
        \noalign{\hrule}
        \lower0pt\hbox\bgroup$
}



\def\blank#1#2{
\hbox to #1{\hfill \vbox to #2{\vfill}}
}

\def\domeq{\trianglerighteq}


\def\strut{\vrule height10pt depth5pt width0pt}

\topmatter
\title Standard modules for tabular algebras
\endtitle

\author R.M. Green \endauthor
\affil 
Department of Mathematics and Statistics\\ Lancaster University\\
Lancaster LA1 4YF\\ England\\
{\it  E-mail:} r.m.green\@lancaster.ac.uk\\
\endaffil

\abstract
We introduce cell modules for the tabular algebras defined in a 
previous work; these modules are analogous to the representations arising
from left Kazhdan--Lusztig cells.  The standard modules of the title are
constructed in an elementary way by suitable tensoring of the cell modules.  
We show how a certain extended affine Hecke algebra of type $A$ equipped 
with its Kazhdan--Lusztig basis is an example of a tabular algebra, and 
verify that in this case our standard modules coincide with other standard 
modules defined in the literature.
\endabstract

\thanks
The author thanks Colorado State University for its hospitality during the
preparation of part of this paper.
\endthanks

\endtopmatter


\centerline{\bf To appear in Algebras and Representation Theory}

\head Introduction \endhead

In their seminal work \cite{{\bf 8}}, Kazhdan and Lusztig showed how left cells
in Hecke algebras may be used to construct representations of the 
algebra.  In favourable cases, such as the Hecke algebras associated to the
symmetric group, all simple representations may be constructed in terms of
these cell representations.  
Graham and Lehrer \cite{{\bf 2}} developed this
idea further by defining cellular algebras in terms of multiplicative 
properties of a basis.  One of the most important results of \cite{{\bf 2}} is 
a classification of the simple modules of a cellular algebra as explicitly
described quotients of the cell modules.

Tabular algebras were introduced by the author in \cite{{\bf 4}} as a class of
associative $\zed[v, v^{-1}]$-algebras equipped with distinguished bases
(tabular bases) 
and satisfying certain axioms.  Although tabular algebras are defined in
terms of a seemingly complicated table datum, the main results of \cite{{\bf 6}} 
show that, under mild assumptions, this table datum may be recovered up to 
isomorphism (in a sense made precise in \cite{{\bf 6}}) from the distinguished 
basis.  There are many natural examples of tabular bases given in \cite{{\bf 4}};
these include the natural basis for Brauer's centralizer algebra and
various bases arising from Kazhdan--Lusztig type constructions (IC bases).

In this paper, we develop the analogue of cell representations for tabular 
algebras.  The axioms for a tabular algebra mean that these are very easy 
to define, and that a cell module inherits a nice basis from the 
corresponding algebra.  By suitable tensoring of these cell modules, we
define in \S1.3 the standard modules for a tabular algebra.  These are
constructed in a simple way from the table datum, but their properties may
be subtle as we shall explain in \S4.

The left, right and two-sided cells of a tabular algebra may be described 
purely in terms of the structure constants, just as in the theory of
Kazhdan--Lusztig cells.  Often a left cell and a right cell will intersect in
more than one element; this is not allowed to happen in the case of cellular
algebras, but it causes no problems in our construction.  Another potential 
advantage of tabular algebras over cellular algebras is that it is possible 
to treat some infinite dimensional examples successfully; one of the main 
results of this paper is Theorem 3.4.5, which shows that a certain
extended affine Hecke algebra of type $A$ is an example of a tabular algebra,
and that the Kazhdan--Lusztig basis in this case serves as a tabular basis.
(This is closely related to N. Xi's results in \cite{{\bf 17}}, as we explain 
in \S4.)
Using an equivalent construction of standard modules developed in \S2, we
show in \S4.2 that the standard modules of the aforementioned Hecke algebra
in the sense of tabular algebras agree with the geometrically defined
standard modules in the sense of \cite{{\bf 12}}.

Some of these ideas have been used implicitly by Graham and Lehrer in the 
special case of their construction of two-step nilpotent representations 
of the extended affine Hecke algebra of type $A$ \cite{{\bf 3}}.  They consider a
certain quotient of the extended affine Hecke algebra
and classify the simple modules for this quotient
in terms of what we call the standard modules.  (Both
the affine Hecke algebra in question and the quotient considered in \cite{{\bf 3}}
are tabular algebras.)  There are also other special cases
of this construction in the literature.  This suggests some directions
for further research; we mention these in the closing remarks.

\head 1. Cell modules and standard modules for tabular algebras \endhead

We begin in \S1 by recalling the definition of a tabular algebra from 
\cite{{\bf 4}}, and introducing the general concept of standard modules for 
tabular algebras.

\subhead 1.1 Tabular algebras \endsubhead

Tabular algebras will be constructed from the normalized table algebras 
defined below.

\definition{Definition 1.1.1}
A normalized table algebra is a pair $(\Gamma, B)$, where $\Gamma$ is an 
associative unital $R$-algebra for some $\zed \leq R \leq \complex$ 
and $B = \{b_i : i \in I\}$ is a distinguished basis for
$\Gamma$ such that $1 \in B$, satisfying the following three axioms:

\item{(T1)}{The structure constants of $\Gamma$ with respect to the basis
$B$ lie in $\real^+$, the nonnegative real numbers.}
\item{(T2)}{There is an algebra anti-automorphism $\bar{\ }$ of $\Gamma$ whose
square is the identity and that has the property 
that $b_i \in B \Rightarrow \overline{b_i} \in
B$.  (We define $\overline{i}$ by the condition $\overline{b_i} =
b_{\bar{i}}$.)}
\item{(T3)}{Let $\k(b_i, a)$ be the coefficient of $b_i$ in $a \in \Gamma$.
Then $$
\k(b_m, b_i b_j) = \k(b_i, b_m \overline{b_j})
,$$ for all $i, j, m$.}
\enddefinition

Notice that the table algebra anti-automorphism is determined by the structure
constants: $b_{\bar{i}}$ is the unique basis element with the property that
$1$ occurs with nonzero coefficient in $b_i b_{\bar{i}}$.

Normalized table algebras as defined above are similar to the table algebras
of Arad--Blau \cite{{\bf 1}} and Sunder's discrete hypergroups \cite{{\bf 16}}. Further
details may be found in \cite{{\bf 4}, \S1.1}.

The definition of $\afn$-function below is due to Lusztig.

\definition{Definition 1.1.2}
Let $\A$ be the ring of Laurent polynomials $\zed[v, v^{-1}]$, let
$A$ be an $\A$-algebra and let $\BB$ be an $\A$-basis of $\A$.
For $X, Y, Z \in \BB$, we define the structure constants $g_{X, Y, Z} \in
\A$ by the formula $$
X Y = \sum_Z g_{X, Y, Z} Z
.$$  The $\afn$-function is defined by $$
\afn(Z) = 
\max_{X, Y \in \BB} \deg(g_{X, Y, Z})
,$$ where the degree of a Laurent polynomial is taken to be 
the highest power of $v$ occurring with nonzero coefficient.  We
define $\g_{X, Y, Z} \in \zed$ to be the coefficient of $v^{\afn(Z)}$ in $g_{X,
Y, Z}$; this will be zero if the bound is not achieved.
\enddefinition

We can now give the definition of a tabular algebra.

\definition{Definition 1.1.3}
A {\it tabular algebra} is an
$\A$-algebra $A$, together with a table datum 
$$(\Lambda, \Gamma, B, M, C, *)$$ satisfying axioms (A1)--(A5) below.

\item{(A1)}
{$\Lambda$ is a finite poset.  For each $\l \in \Lambda$, 
$(\Gamma(\l), B(\l))$ is
a normalized table algebra over $\zed$ and
$M(\l)$ is a finite set.  The map $$
C : \coprod_{\l \in \Lambda} \left( M(\l) \times B(\l) \times M(\l)
\right) \rightarrow A
$$ is injective with image an $\A$-basis of $A$.  We assume
that $\Im(C)$ contains a set of mutually orthogonal idempotents 
$\{1_\e : \e \in {\Cal E}\}$ such that 
$A = \sum_{\e, \e' \in {\Cal E}} (1_\e A 1_{\e'})$ and such that for each
$X \in \Im(C)$, we have $X = 1_\e X 1_{\e'}$ for some $\e, \e' \in
{\Cal E}$.
(Typically, the above set of idempotents contains only the identity element 
of $A$.)
A basis arising in
this way is called a {\it tabular basis}.  
}
\item{(A2)}
{If $\l \in \Lambda$, $S, T \in M(\l)$ and $b \in B(\l)$, we write
$C(S, b, T) = C_{S, T}^{b} \in A$.  
Then $*$ is an $\A$-linear involutory anti-automorphism 
of $A$ such that
$(C_{S, T}^{b})^* = C_{T, S}^{\overline{b}}$, where $\bar{\ }$ is the
table algebra anti-automorphism of $(\Gamma(\l), B(\l))$.
If $g \in \complex(v) \otimes_\zed \Gamma(\l)$ is such that 
$g = \sum_{b_i \in B(\l)} c_i b_i$ for some scalars $c_i$ 
(possibly involving $v$), we write
$C_{S, T}^g \in \complex(v)\otimes_\A A$ 
as shorthand for $\sum_{b_i \in B(\l)} c_i C_{S, T}^{b_i}$.  We write
$\bc_\l$ for the image under $C$ of $M(\l) \times B(\l) \times M(\l)$; it 
turns out \cite{{\bf 4}, Proposition 2.3.1} that the $\afn$-function is 
constant on each set $\bc_\l$.}
\item{(A3)}
{If $\l \in \Lambda$, $g \in \Gamma(\l)$ and $S, T \in M(\l)$ then for all 
$a \in A$ we have $$
a . C_{S, T}^{g} \equiv \sum_{S' \in M(\l)} C_{S', T}^{r_a(S', S) g}
\mod A(< \l),
$$ where  $r_a (S', S) \in \Gamma(\l)[v, v^{-1}] = \A \otimes_\zed
\Gamma(\l)$ is independent of $T$ and of $g$ and $A(< \l)$ is the
$\A$-submodule of $A$ generated by the set $\bigcup_{\mu < \l} \bc_\mu$.}
\item{(A4)}{Let $K = C_{S, T}^b$, $K' = C_{U, V}^{b'}$ and 
$K'' = C_{X, Y}^{b''}$ lie in $\Im(C)$.  Then the
maximum bound for $\deg(g_{K, K', K''})$
in Definition 1.1.2 is achieved if and only if $X = S$, $T = U$, $Y =
V$ and $\kappa(b'', bb') \ne 0$ (where $\kappa$ is as in axiom (T3)).  
If these conditions all hold and
furthermore $b = b' = b'' = 1$, we require $\g_{K, K', K''} = 1$.}
\item{(A5)}{There exists an $\A$-linear function $\t : A \ra \A$
(the {\it tabular trace}), such that $\t(x) = \t(x^*)$ for all $x \in
A$ and $\t(xy) = \t(yx)$ for all $x, y \in A$, that has the 
property that for every
$\l \in \Lambda$, $S, T \in M(\l)$, $b \in B(\l)$ and $X = C_{S,
T}^b$, we have $$
\t(v^{\afn(X)} X) = 
\cases 1 \mod v^{-1} \A^- & \text{ if } S = T \text{ and } b = 1,\cr
0 \mod v^{-1} \A^- & \text{ otherwise.} \cr
\endcases
$$  Here, $\A^- := \zed[v^{-1}]$.  We call the elements $C_{S, S}^1$
{\it distinguished involutions}.}
\enddefinition

\remark{Remark 1.1.4}
In \cite{{\bf 4}}, a tabular algebra is only required to satisfy axioms (A1)--(A3),
and an algebra satisfying all five axioms is called a ``tabular algebra with 
trace''.  However, all the tabular algebras in this paper (and all the most
interesting examples) are tabular
algebras with trace, so we use the term ``tabular algebra'' with this
narrower meaning.
\endremark

\remark{Remark 1.1.5}
The table datum for a tabular algebra $A$ with tabular 
basis $\BB$ may be reconstructed
up to isomorphism (in a sense made precise in \cite{{\bf 6}}) from the structure
constants.  This means that there is no loss in considering tabular algebras
to be pairs $(A, \BB)$; we shall look at this in more detail in \S3 in the case
where $A$ is an affine Hecke algebra and $\BB$ is its Kazhdan--Lusztig basis.
\endremark

\subhead 1.2 Cell modules \endsubhead

We now introduce cell modules for tabular algebras, the idea of which is
implicit in the formulation
of axiom (A3).  Cell modules are the analogues of the left cell representations
in \cite{{\bf 8}} and of the cell modules in \cite{{\bf 2}, Definition 2.1}.

\definition{Definition 1.2.1}
Let $A$ be a tabular algebra with table datum $(\Lambda, \Gamma, B, M, C, *)$.
For each $\l \in \Lambda$, we define the (left) $A$-module $W(\l)$ as follows:
$W(\l)$ is a free $\A$-module with basis $\{C_S^g : S \in M(\l), g \in B(\l)\}$
and $A$-action defined by $$
a C_S^g = \sum_{S' \in M(\l)} C_{S'}^{r_a(S', S)g}
$$ for all $a \in A$.  Here, $r_a(S', S)$ is as in axiom 
(A3), and the notational shorthand is analogous to that in axiom (A3).  
This is called a (left) cell module for $A$, and the representation it
affords is called a (left) cell representation.
\enddefinition

\remark{Remark 1.2.2}
It is clear from axiom (A3) that this module action is well defined.  In fact,
more is true: the cell modules may be recovered (as modules with distinguished
bases) from the tabular basis $\BB$ (see Corollary 1.2.6 below).

We could also define right cell modules by applying the tabular 
anti-auto\-morph\-ism, $*$.
\endremark

We now define the left cells of a tabular algebra; this is the one-sided 
version of the construction in \cite{{\bf 4}, \S3.1} referred to in \cite{{\bf 4},
Remark 3.1.2}.

\definition{Definition 1.2.3}
Let $A$ be a tabular algebra.  We introduce a relation, $\preceq_L$, on
the tabular basis by stipulating that $X' \preceq_L X$ if $X'$ appears
with nonzero coefficient in $KX$ for some tabular basis
element $K$.
\enddefinition

The proof of Lemma 1.2.5 will need the following

\definition{Definition 1.2.4}
Let $A$ be a tabular algebra with table datum $(\Lambda, \Gamma, B, M,
C, *)$.  Let $\l \in \Lambda$ and $S, T, U, V \in M(\l)$.  We define
$\lan T, U \ran \in \Gamma(\l)[v, v^{-1}]$ by the condition $$
C_{S, T}^1 C_{U, V}^{1} \equiv C_{S, V}^{\lan T, U \ran} \mod A(< \l)
.$$  If $b \in B(\l)$, we define $\lan T, U \ran_b \in \A$ to be the 
coefficient of $b$ in $\lan T, U \ran$.
\enddefinition

Definition 1.2.4 is sound by \cite{{\bf 4}, Lemma 2.2.2}.

\proclaim{Lemma 1.2.5}
Let $A$ be a tabular algebra with table datum 
$(\Lambda, \Gamma, B, M, C, *)$.  Let $\preceq_L^t$ be the transitive
extension of the relation $\preceq_L$ of Definition 1.2.3.  The
relation $\sim_L$ on $\Im(C)$ defined by $Y \sim_L Z$ if and
only if $Y \preceq_L^t Z$ and $Z \preceq_L^t Y$ is an equivalence
relation.  The equivalence classes are known as left cells.  Two basis elements
$C_{T, U}^b$ and $C_{V, W}^{b'}$ are in the same left cell if and only if
$U = W$ as elements of $\coprod_{\l \in \Lambda} M(\l)$.
\endproclaim

\demo{Proof}
The proof is a simple adaptation of the proof of \cite{{\bf 4}, Proposition 3.1.3}.

The idempotent condition in axiom (A1) shows that $\sim_L$ is reflexive.
Since $\sim_L$ is clearly symmetric and transitive, it is
an equivalence relation.  

Let $Y = C_{T, U}^b$ and $Z = C_{V, U}^{b'}$ be basis elements such that
$T, U, V \in M(\l)$ for the same $\l \in \Lambda$; we will show
that $Y \sim_L Z$.  Now $$
Y^* Y = C_{U, T}^{\overline{b}} C_{T, U}^b \equiv 
C_{U, U}^{\overline{b} \lan T, T \ran b} \mod A(< \l)
,$$ which, by \cite{{\bf 4}, Lemma 2.2.3 (i)}, contains $C_{U, U}^1$ with nonzero
coefficient of degree $\afn(\l)$, where $\afn(\l)$ is the $\afn$-value of
any basis element in $\bc_\l$.  (The function $\afn$ is well defined on 
$\Lambda$ by \cite{{\bf 4}, Proposition 2.3.1}, and the sets $\bc_\l$ are as 
given in axiom (A2).)
There is a similar converse
statement: $C_{T, U}^b C_{U, U}^1$ contains $Y$ with nonzero
coefficient.  This shows that $Y \sim_L C_{U, U}^1$.  Similarly, we have
$Z \sim_L C_{U, U}^1$ and hence $Y \sim_L Z$ as claimed.

Axiom (A3) shows that if 
$C_{T, U}^b \preceq_L^t C_{V, W}^{b'}$, then either $U = W$ or the two
basis elements are in different two-sided cells $\bc_\l$.
The partial order on $\Lambda$ means that elements in different two-sided
cells must be in different left cells.  It follows that the equivalence classes
are as claimed.
\qed\enddemo

\proclaim{Corollary 1.2.6}
The cell modules (considered as left $A$-modules with distinguished
bases) are determined by the tabular basis.
\endproclaim

\demo{Proof}
By Definition 1.2.3 and Lemma 1.2.5, the left cells depend only on the tabular
basis $\BB$ and not on the rest of the table datum.  Given an left cell 
$L \subset \BB$, we
may construct a module with basis $\{l + A(< \l) : l \in L\}$ with the obvious
left $A$-action.  This is well defined by the relationship between left cells
and two-sided cells mentioned 
in the last paragraph of the proof of Lemma 1.2.5.  

Any given table datum will
provide an isomorphism of based $\A$-modules between this module and the
module $W(\l)$ of Definition 1.2.1: an element $C_{S, T}^g$ in $L$ corresponds
to the basis element $C_S^g$.
\qed\enddemo

\remark{Remark 1.2.7}
The labelling of the basis elements of $W(\l)$ is highly dependent on the table
datum: for example, there is generally no way to identify which basis elements
of $W(\l)$ are of the form $C_S^1$.
We do not pursue this, but it may be proved using the results of \cite{{\bf 6}}.
\endremark

\subhead 1.3 Standard modules \endsubhead

One of the important properties of the cell modules $W(\l)$ is
that, as well as being left modules for the tabular algebra $A$, they are also
right modules for the hypergroup $\Gamma(\l)$.  We now show that these actions
commute with each other.  

The above notation will be fixed throughout \S1.3, as will the table 
datum for $A$.

\proclaim{Lemma 1.3.1}
For any $g \in \Gamma(\l)$, the $\A$-linear map $\psi_g : W(\l) \ra W(\l)$
defined by $$\psi_g(C_S^b) := C_S^{bg}$$ is a homomorphism of left $A$-modules.
This gives $W(\l)$ the structure of a free right $\Gamma(\l)$-module, and the
$A$-action and $\Gamma(\l)$-action on $W(\l)$ commute with each other.
\endproclaim

\demo{Proof}
The map $\psi_g$ is a homomorphism by axiom (A3), and the independence of
$r_a(S', S)$ from $g$ in that axiom shows that the two actions commute as
claimed.  Freeness follows from the fact that $C_S^1 . g = C_S^g$.
\qed\enddemo

\remark{Remark 1.3.2}
The $\Gamma(\l)$-action defined in Lemma 1.3.1 does in general depend on the
table datum chosen.  (Contrast this with Corollary 1.2.6.)
\endremark

Although Lemma 1.3.1 is valid integrally (\idest over $\A$), it will
be necessary for some constructions to take the base ring to be a field.

\definition{Definition 1.3.3}
Let $k$ be a field and let $r \in k^*$.  Then $k$ is naturally a unital 
$\A$-module, with $v$ acting as multiplication by $r$.  If $W(\l)$ is a 
cell module for a tabular algebra $A$ with table datum $(\Lambda, \Gamma, B,
M, C, *)$, we write $W(\l)_{(k, r)}$ for $k \otimes_\A W(\l)$ and
$A_{(k, r)}$ for $k \otimes_\A A$, where the $\A$-module structure of $k$ is
as above.  It is clear that $W(\l)_{(k, r)}$ is a left $A_{(k, r)}$-module.
Similarly, we write $\Gamma(\l)_{(k, r)}$ for $k \otimes_\A \A \otimes_\zed 
\Gamma(\l)$.  Just as in Lemma 1.3.1, $W(\l)_{(k, r)}$ is an 
$A_{(k, r)}$--$\Gamma(\l)_{(k, r)}$-bimodule in which the two actions commute.
\enddefinition

Standard modules are defined as follows.

\definition{Definition 1.3.4}
Maintain the notation of Definition 1.3.3.  A {\it standard module} for the
tabular algebra $A$ is a left $A_{(k, r)}$-module
$W(\l)_{(k, r)} \otimes_{\Gamma(\l)_{(k, r)}} N$, where $N$ is a simple 
$\Gamma(\l)_{(k, r)}$-module.
\enddefinition

\remark{Remark 1.3.5}
The labelling of the standard $A$-modules by pairs $(\l, N)$ may depend 
on the table datum because of Remark 1.3.2, but we shall see that 
the set of isomorphism classes of
standard modules is determined by the tabular basis (Theorem 2.2.3 (iv)).
\endremark

\head 2. Another construction of standard modules \endhead

In order to develop the properties of standard modules, we give a second
construction in terms of the asymptotic tabular algebra.  This construction
will be valuable in \S3 when we study the affine Hecke algebra of $GL_n$ as a
tabular algebra.

\vfill\eject
\subhead 2.1 Asymptotic tabular algebras and bimodules \endsubhead

The asymptotic versions of tabular algebras were introduced in 
\cite{{\bf 4}, \S3}, using methods from \cite{{\bf 13}}.

\definition{Definition 2.1.1}
Let $A$ be a tabular algebra with trace, and maintain the usual notation.
Define $\widehat{X} := v^{-\afn(X)} X$ for any
tabular basis element $X \in \Im(C)$.  
The free $\A^-$-module $A_\l^-$ is
defined to be generated by the elements $\{\widehat{X} : X \in
\bc_\l\}$.  We set $t_X$ to be the image of $\widehat{X}$ in $$
A_\l = A_\l^\infty := {{A_\l^-} \over {v^{-1} A_\l^-}}
.$$  The latter is a $\zed$-algebra with basis $\{t_X : X \in
\bc_\l\}$ and structure constants $$
t_X t_{X'} = \sum_{X'' \in \bc} \g_{X, X', X''} t_{X''}
,$$ where the $\g_{X, X', X''} \in \zed$ are as in Definition 1.1.2.
We also set $$
A^\infty := \bigoplus_{\l \in \Lambda} A_\l^\infty
;$$ this is a $\zed$-algebra with basis $\{t_X : X \in \Im(C)\}$ called
the {\it asymptotic tabular algebra}.  We will use the notation 
$A_{\l (k, r)}^\infty$ and $A_{(k, r)}^\infty$ in the usual way (see
Definition 1.3.3) to denote change of scalars by tensoring with
$k \otimes_\A \A \otimes_\zed -$.
\enddefinition
 
We now introduce a certain bimodule that will be helpful in our second
construction of standard modules.

\definition{Definition 2.1.2}
Let $A$ be a tabular algebra (over $\A$) with table datum 
\newline $(\Lambda, \Gamma, B, M, C, *)$.  We define $E(\l)$ to be the 
free $\A$-module with basis $\bc_\l$ and we give $E(\l)$ the structure of
an $A$--$A_\l^\infty$-bimodule as follows.

The left $A$-module structure of $E(\l)$ is the natural one arising from
the identification of $E(\l)$ with $A(\leq \l)/A(< \l)$ given by
sending a basis element $X \in E(\l)$ to $X + A(< \l)$.
Note that $E(\l)$ is isomorphic as an $A$-module to the direct sum of 
$|M(\l)|$ copies of $W(\l)$.  
On the other hand, identifying basis elements of $E(\l)$ with basis 
elements of $A_\l^\infty$ by the correspondence $X \leftrightarrow t_X$ 
defines the right $A_\l^\infty$-action on $E(\l)$ via the regular 
representation over $\A$.
\enddefinition

\proclaim{Lemma 2.1.3}
The actions of $A$ and $A_\l^\infty$ on $E(\l)$ commute with each other.
\endproclaim

\demo{Proof}
Recall from \cite{{\bf 4}, Theorem 3.2.4 (i)} that 
$A_\l^\infty \cong M_{|M(\l)|}(\Gamma(\l))$ as $\zed$-algebras; the 
isomorphism is provided by the table datum.  It is thus enough
to check that the left $A$-action on $E(\l)$ commutes with the action
of the element of $A_\l^\infty$ corresponding to $e_{ST} \otimes_\zed b$ for 
$S, T \in M(\l)$ and $b \in B(\l)$.  The latter action is defined by $$
C_{U, V}^{b'} . (e_{ST} \otimes_\zed b) = \delta_{VS} C_{U, T}^{b'b}
,$$ and axiom (A3) is precisely what is needed for the conclusion to hold.
\qed\enddemo

Since $A_\l^\infty$ is a ring with identity ($\Gamma(\l)$ has
identity and $A_\l^\infty$ is a matrix ring over it) and $E(\l)$ affords the
regular representation of $A_\l^\infty$, we have 
$\End_{A_\l^\infty}(E_\l) \cong A_\l^\infty$,
where the isomorphism is given by left multiplication.
Lemma 2.1.3 thus induces a homomorphism $\Phi_\l : A \ra \A 
\otimes_\zed A_\l^\infty$ as follows.

\proclaim{Lemma 2.1.4}
The homomorphism $\Phi_\l : A \ra \A \otimes_\zed A_\l^\infty$ induced by 
Lemma 2.1.3
is given by $$\Phi_\l(X) = \sum_{D \in \D_\l, Z \in \bc_\l} g_{X, D, Z} t_Z
,$$ where $\D_\l := \{ C_{S, S}^1 : S \in M(\l) \}$.
\endproclaim

\demo{Proof}
We identify $E_\l$ with $\A \otimes_Z A_\l^\infty$ as in Definition 2.1.2.  
Let $X$ be any basis element of $\Im(C)$ (we do not assume that $X 
\in \bc_\l$).  An element $\phi \in \End_{A_\l^\infty}(A_\l^\infty)$ is 
defined by the value of $\phi(1)$,
so the Lemma will follow once we calculate the value of $X(1)$.  Recall
that the identity element of $A_\l^\infty$ is $\sum_{D \in \D(\l)} t_D$; 
this is Lusztig's property $P_2$ in \cite{{\bf 4}, \S3.2}.

Using the left action of $A$ on $E(\l)$, we calculate that $$
X . \left( \sum_{D \in \D_\l} D \right) =  
\sum_{D \in \D_\l, Z \in \bc_\l} g_{X, D, Z} Z
,$$ where both sides of the equation are considered as elements of $E(\l)$.
Considering $E(\l)$ as the regular representation of $A_\l^\infty$ again, 
we see
that the endomorphism of $E(\l)$ induced by $X$ can also be obtained by
left multiplication by the element $$
\sum_{D \in \D_\l, Z \in \bc_\l} g_{X, D, Z} t_Z
$$ of $A_\l^\infty$.  The conclusion follows.
\qed\enddemo

The homomorphism $\Phi_\l$ first appeared in the work of Lusztig; see
\cite{{\bf 13}, \S1.8} and \cite{{\bf 4}, \S3.2}.  It will be useful to have an
alternative formula for $\Phi_\l$ using the notation $\lan T, U \ran_b$ given
in Definition 1.2.4.  The translation of Lemma 2.1.4 into this notation is
as follows.

\proclaim{Lemma 2.1.5}
The homomorphism $\Phi_\l : A \ra \A \otimes_\zed A_\l^\infty$ 
is given by $$\Phi_\l(a) = \sum_{S, S' \in M(\l), \ b \in B(\l)} 
(r_a(S', S))_b \  t_{S', S}^b
.$$  Here, $t_{S', S}^b = e_{S', S} \otimes b$ is the element $t_Z$ 
where $Z = C_{S', S}^b$. \qed
\endproclaim

\subhead 2.2 Constructing $A$-modules using the asymptotic tabular
algebra \endsubhead

We now show how the cell modules of \S1.2 may be defined starting from the
asymptotic tabular algebra.

\definition{Definition 2.2.1}
Let $A$ be a tabular algebra (over $\A$) with table datum 
\newline $(\Lambda, \Gamma, B, M, C, *)$.  For each $\l \in \Lambda$, we
define $V(\l)$ to be the $\zed$-module with basis $\{v_S : S \in M(\l)\}$.
The $\zed$-module $W'(\l) := V(\l) \otimes_\zed \Gamma(\l)$ has the 
structure of a left $A_\l^\infty$-module via $$
(e_{ST} \otimes b) . (v_U \otimes g) = \d_{TU} v_S \otimes bg
,$$ using the isomorphism $A_\l^\infty \cong M_{|M(\l)|}(\Gamma(\l))$ from 
\cite{{\bf 4}, Theorem 3.2.4 (i)}.  The $\A$-module $\A \otimes_\zed W'(\l)$ 
becomes an $A$-module using the homomorphism $\Phi_\l$.  If $k$ is a field
and $r \in k^*$, we write $W'(\l)_{(k, r)}$ for the $k$-module 
$k \otimes_\A \A \otimes_\zed W'(\l)$; this is a left $A_{(k, r)}$-module
in the obvious way.
\enddefinition

\proclaim{Proposition 2.2.2}
Let $A$ be a tabular algebra (over $\A$) with table datum
\newline $(\Lambda, \Gamma, B, M, C, *)$, and let $\l \in \Lambda$.
Let $\th : W(\l) \ra \A \otimes_\zed W'(\l)$ be the $\A$-module isomorphism 
sending $C_T^b$ to $1 \otimes v_T \otimes b$.  Then $\th$ is an isomorphism of
$A$--$\Gamma(\l)$-bimodules.
\endproclaim

\demo{Proof}
Let us first note that the actions of $A$ and $\Gamma(\l)$ on $\A \otimes_\zed
W'(\l)$ commute.  This is because, as one can check, the action of 
$A_\l^\infty$ on 
$\A \otimes_\zed W'(\l)$ commutes with the action of $\Gamma(\l)$.  The same
is true for $W(\l)$ by Lemma 2.1.3.

It is clear that $\th$ is compatible with the $\Gamma(\l)$ actions on both
modules, because $\th(C_T^{bg}) = v_T \otimes bg$ by linearity.

Let $a \in A$.  Then using Lemma 2.1.5, we have $$\eqalign{
a . (v_T \otimes g) &= \Phi_\l(a) . (v_T \otimes g) \cr
&= \left( \sum_{S, S' \in M(\l), \ b \in B(\l)} 
(r_a(S', S))_b \  t_{S', S}^b \right) .  (v_T \otimes g) \cr
&= \sum_{S' \in M(\l), \ b \in B(\l)} 
(r_a(S', T))_b \  v_{S'} \otimes bg.\cr
}$$  On the other hand, $
a . C_T^g = \sum_{S' \in M(\l)} C_{S'}^{r_a(S', T)g}
$ from the definition of $W(\l)$.  Because $r_a(S', T) = \sum_{b \in B(\l)}
r_a(S', T)_b \ b$, we have $$
a . C_T^g = \sum_{S' \in M(\l)} (r_a(S', T))_b \ C_{S'}^{bg}
.$$  It now follows that $\th$ is a homomorphism of left $A$-modules.
\qed\enddemo

\proclaim{Theorem 2.2.3}
Let $k$ be a field and let $r \in k^*$.  
Let $A$ be a tabular algebra with table datum 
$(\Lambda, \Gamma, B, M, C, *)$, and let $\l \in \Lambda$.
Let $$\th_{(k, r)} : W(\l)_{(k, r)} \ra W'(\l)_{(k, r)}$$ be the 
$\A$-module isomorphism 
sending $C_S^b$ to $v_S \otimes b$.  Then
\item{\rm (i)}{the map $\th_{(k, r)}$ is an isomorphism of
$A_{(k, r)}$--$\Gamma(\l)_{(k, r)}$-bimodules;}
\item{\rm (ii)}{every standard module for $A_{(k, r)}$ is isomorphic to 
$W'(\l)_{(k, r)} \otimes_{\Gamma(\l)_{(k, r)}} N$ for some simple
$\Gamma(\l)_{(k, r)}$-module $N$;}
\item{\rm (iii)}{every standard module for $A_{(k, r)}$ is isomorphic 
to a simple $A_{(k, r)}^\infty$-module $N'$ regarded as an $A_{(k, r)}$-module
by setting $a . n = \Phi_\l(a) . n$ for some $\l \in \Lambda$;}
\item{\rm (iv)}{the isomorphism classes of standard modules for 
$A_{(k, r)}$ are determined by the tabular basis.}
\endproclaim

\demo{Proof}
Part (i) follows from Proposition 2.2.2 after a suitable change of scalars,
and part (ii) follows from (i) and Definition 1.3.4.  

For (iii), note that
the algebras $A_{\l (k, r)}$ and $\Gamma(\l)_{(k, r)}$ are Morita equivalent,
since the first is isomorphic to an $|M(\l)|$ by $|M(\l)|$ 
matrix ring over the second.  This means
that the functor $V(\l)_{(k, r)} \otimes_k \Gamma(\l)_{(k, r)} 
\otimes_{\Gamma(\l)_{(k, r)}} -$ identifies simple 
$\Gamma(\l)_{(k, r)}$-modules with 
simple $A_{\l (k, r)}$-modules.  The conclusion now follows from (ii)
once we recall that $A_{\l (k, r)}$ is a direct summand of 
$A_{(k, r)}^\infty$.

The definition of standard modules given in (iii) is determined by the
tabular basis.  At first sight it seems that
the definition depends on $\Lambda$, but in fact it only depends on the
two-sided cells in $A$, the corresponding homomorphisms $\Phi_\l$ and
the asymptotic algebras $A_\l^\infty$.  The maps $\Phi_\l$ of Lemma 2.1.4
depend only on the tabular basis since the set ${\Cal D}_\l$ depends only
on the basis (by axiom (A5)).  The two-sided cells and asymptotic algebras 
are determined by the tabular basis by \cite{{\bf 6}, Proposition 2.3.3,
Corollary 2.3.4}.  These observations prove (iv).
\qed\enddemo

\head 3. The affine Hecke algebra of $GL_n$ \endhead

In \S3, we show that the asymptotic Hecke algebra of $GL_n$ equipped with
the Kazhdan--Lusztig basis is a tabular algebra, thus answering in the
affirmative a question raised in \cite{{\bf 4}, \S7}.  The proof of this result
relies on Xi's solution \cite{{\bf 17}} to Lusztig's conjecture \cite{{\bf 12}, Conjecture
10.5} on the structure of the associated asymptotic Hecke algebra.

\subhead 3.1 Hecke algebras \endsubhead

We will define Hecke algebras for extended Coxeter groups 
following the notation of \cite{{\bf 17}, \S1.1}.

\definition{Definition 3.1.1}
Let $(W',S)$ be a Coxeter system with $S$ the set of simple
reflections. Assume that a commutative group $\Omega$ acts on
$(W',S)$, and consider the {\it extended Coxeter group}
$W=\Omega\ltimes W'$. The length function $\ell$ on $W'$ and the
Bruhat--Chevalley order $\leq$ on $W'$ are extended to $W$ by stipulating
that $\ell(\omega w):=\ell(w)$, and $\omega w\leq
\omega' u$ if and only if $\omega=\omega'$ and $w\leq u$, where
$\omega,\omega'$ are in $\Omega$ and $w,u$ are in $W'$.

The Hecke algebra $\H$ of $(W,S)$
over $\A$ with parameter $v^2$ is an associative
algebra over $\A$, with free $\A$-basis $\{T_w : w\in W\}$ and defining
relations $$\eqalign{
(T_s-v^2)(T_s+1) &=0 \  \text{ if } s\in S; \cr
T_wT_u &= T_{wu} \ \text{ if } \ell(wu)=\ell(w)+\ell(u).
}$$

Let $\bar{ } : \A \ra \A$ be the $\zed$-linear ring homomorphism 
$\A$ defined by $\bar
v=v^{-1}$. Then we have a bar involution of $\H$ defined by $$
\overline{\sum a_wT_w}=\sum \bar a_wT^{-1}_{w^{-1}},\qquad
a_w\in\A.$$ For each $w\in W$ there is a unique element $C_w$ in
$\H$ such that $\overline{C_w}=C_w$ and $$
C_w=v^{-\ell(w)}\sum_{y\leq w}P_{y,w}(v^2)T_y
,$$ where $P_{y,w}$ is a polynomial in $v$ of
degree $\leq \frac 12(\ell(w)-\ell(y)-1)$ if $\ell(w)>\ell(y)$ and
$P_{w,w}=1$.

The basis $\{C_w : w\in W\}$ is called the Kazhdan--Lusztig basis of the 
Hecke algebra $H$ and the polynomials $P_{y,w}$ are
the Kazhdan--Lusztig polynomials.
\enddefinition

We can now introduce the affine Hecke algebra of $GL_n$.  Xi calls this the
``extended affine Hecke algebra associated to $SL_n(\complex)$''
\cite{{\bf 17}, \S8.4}, and Graham--Lehrer \cite{{\bf 3}} 
call it the ``extended affine Hecke algebra of type $A$''.

\definition{Definition 3.1.2}
The affine Hecke algebra of $GL_n$ arises from the construction in Definition
3.1.1 by setting $(W', S)$ to be the Coxeter system of type $\widehat A_{n-1}$
(for $n \geq 3$) and $\Omega$ to be the cyclic group $\zed_n$ acting by
rotations of the Coxeter graph.
\enddefinition

\remark{Remark 3.1.3}
The Kazhdan--Lusztig basis arising from the extended Coxeter group $(W, S)$ 
is closely related to that arising from the Coxeter group $(W', S)$.  The
basis elements of the larger Hecke algebra are precisely those of the form 
$T_\omega C_u$, where $\omega \in \Omega$ and $u \in W'$; this follows from
the characterization of the basis given in Definition 3.1.1 because $T_\omega$
is invariant under the involution.  In particular, $C_\omega = T_\omega$.
It also follows easily from Definition 3.1.1 that 
$T_\omega C_u T_{\omega^{-1}}$ is a
Kazhdan--Lusztig basis element, namely $C_{\omega u \omega^{-1}}$ (consider
the leading term $T_u$ of $C_u$).  These properties of the basis are well 
known; see for example the proof of \cite{{\bf 17}, Proposition 1.4.6 (b)}.
\endremark

\remark{Remark 3.1.4}
We note that there are other quite different presentations of the algebra
given in Definition 3.1.2, but the set-up given is 
convenient for studying Kazhdan--Lusztig theory.  
\endremark

The following standard result will be useful later.

\proclaim{Lemma 3.1.5}
Let $(W, S)$ be an extended Coxeter group.  
The $\A$ linear map $*$ on $\H$ that sends $T_w$ to $T_{w^{-1}}$ is an algebra
anti-automorphism, and we have $C_w^* = C_{w^{-1}}$ for all $w \in W$.
\endproclaim

\demo{Proof}
The first assertion comes from symmetry properties of the relations for $\H$
given in Definition 3.1.1 and the fact that $\ell(w) = \ell(w^{-1})$ for 
all $w \in W$.  Applying $*$ to the defining expression for $C_w$ in
Definition 3.1.1 shows that $C_w^*$ satisfies the required degree bounds.

It only remains to show that $\overline{C_w^*} = C_w^*$.
It follows from the definition of $\H$ that it is generated as an 
$\A$-algebra by $\{C_s : s \in S\}$ and $\{C_\omega: \omega \in \Omega\}$
and it is enough to prove that $*$ and $\bar{\ }$ commute on $\H$.  This 
follows because $* \circ \bar{\ }$ and $\bar{\ } \circ *$ are both
$\A$-antilinear ring homomorphisms fixing the algebra generating set
just given.
\qed\enddemo

We now turn our attention to a particularly important set of involutions 
in the Coxeter group $W'$.

\definition{Definition 3.1.6}
Let $(W', S)$ be a Coxeter group.  Following Lusztig \cite{{\bf 10}, \S1.3}, we
find that $\afn(z) \leq \ell(z) - 2 \deg P_{e, z}$, and we 
define the set of {\it distinguished involutions} of $W'$ to 
be $$
{\Cal D} = \{ z \in W' : 2 \deg P_{e, z} = \ell(z) - \afn(z) \}
.$$
\enddefinition

It is not immediately clear that the distinguished involutions are involutions,
but this is proved in \cite{{\bf 10}, Proposition 1.4}.  The terminology of 
Definition 3.1.6 will eventually be seen to be compatible with the 
distinguished involutions mentioned in axiom (A5) for a tabular algebra.

\subhead 3.2 Kazhdan--Lusztig cells \endsubhead

For the rest of \S3, we shall
usually denote the affine Hecke algebra of $GL_n$ by $A$, the 
Kazhdan--Lusztig basis by $\BB$ and the associated extended Coxeter group 
by $W$.  We aim to show that $\BB$ is a tabular basis for $A$.
The table datum for $(A, \BB)$ will be given in terms of the two-sided
Kazhdan--Lusztig cells of the Hecke algebra.  These may be defined in terms
similar to Definition 1.2.3.

\definition{Definition 3.2.1}
Let $(W', S)$ be a Coxeter group and let 
$x, w \in W'$.  We write $x \leq_L w$ if there is a chain $$ x
= x_0, x_1, \ldots, x_r = w ,$$ possibly with $r = 0$, such that
for each $i<r$, $C_{x_i}$ occurs with nonzero coefficient in the
linear expansion of $C_s C_{x_{i+1}}$ for some $s \in S$ such
that $s x_{i+1} > x_{i+1}$.

The preorder $\leq_R$ on $W'$ is defined by the
condition $x \leq_R w \Leftrightarrow x^{-1} \leq_L w^{-1}$, and
the preorder $\leq_{LR}$ is that generated by $\leq_L$ and
$\leq_R$.

If $W = \Omega\ltimes W'$ is an extended Coxeter group, $w, u \in W'$ 
and $\omega,\omega'$ in $\Omega$, we say that $\omega w \leq_L \omega'
u$ (respectively, $\omega w \leq_R \omega' u$, $\omega w \leq_{LR} \omega' u$) 
if $w \leq_L u$ (respectively $w \leq_R u,\ w \leq_{LR} u$).

The transitive preorder $\leq_L$ yields an equivalence relation $\sim_L$
on $W$ or $W'$ (where $x \sim_L w$ if and only if $x \leq_L w$ and $w
\leq_L x$) whose equivalence classes are called the {\it left
cells} of $W$.  The preorders $\leq_R$ and $\leq_{LR}$ 
yield equivalence relations $\sim_R$ and
$\sim_{LR}$ on $W$ whose equivalence classes are called {\it right
cells} and {\it two-sided cells}, respectively.  Note that the left 
(respectively right, two-sided) cells are partially ordered by $\leq_L$
(respectively $\leq_R$, $\leq_{LR}$).
\enddefinition

The above definition agrees with the original definition in \cite{{\bf 8}} for
the Coxeter group $W'$; see \cite{{\bf 7}, Remark 1.2.3} for an explanation.

The next result is immediate from the definitions.

\proclaim{Lemma 3.2.2}
Let $(W, S)$ be an extended Coxeter group corresponding to the Coxeter group
$(W', S)$.  The two-sided (respectively left, right) cells of
$W$ are precisely the subsets of the form $\{\omega c : \omega \in \Omega,
c \in \bc\}$, where
$\bc$ is a two-sided (respectively left, right) cell of $W'$.
\qed\endproclaim

The following well known result will be useful in connection with tabular
bases.

\proclaim{Lemma 3.2.3}
Let $(W, S)$ be an extended Coxeter group with Hecke algebra $\H$ as defined
in \S3.1.  For any $\xi,\xi'\in \H$ and $w\in W$ we write $$
\xi C_w \xi'=\sum_{z \in W} h_z C_z
$$ with $h_{z} \in \A$.  Then $h_z = 0$ unless $z \leq_{LR} w$.  
If $\xi' = 1$ (respectively, $\xi = 1$), we have $h_z = 0$ unless $z \leq_L w$
(respectively, $z \leq_R w$).  Furthermore, the preorders $\leq_L$, $\leq_R$
and $\leq_{LR}$ are characterized by these conditions.
\endproclaim

\demo{Proof}
The last assertion follows from Definition 3.2.1, because $C_s C_w = 
(v + v^{-1}) C_w$ if $s, w \in W'$ and $sw < w$ (where $W'$ is the Coxeter
group corresponding to $W$).
The other assertions were proved in \cite{{\bf 9}, \S4.3} in the case
where $W$ is a Coxeter group.
For the general case, we need only check the case where the elements $\xi$
and $\xi'$ are of the form $T_\omega$ for $\omega \in \Omega$, and in this
case the result is immediate from Definition 3.2.1.
\qed\enddemo

For the rest of this section, let $(W', S)$ be a Coxeter group of
type $\widehat{A}_{n-1}$ and $(W, S)$ be the extended Coxeter group described
in Definition 3.1.2.  

In this case, the cells have a nice description 
established by Shi \cite{{\bf 14}}, who showed that the two-sided cells are labelled
by the partitions of $n$.  This description is given by a 
map, $\s$, from $W'$ to the set of partitions of $n$, whose fibres are the
two-sided cells.  (The reader is referred to \cite{{\bf 14}} for the definition 
of $\s$, which we do not require here.)  We extend $\s$ to a map from 
$W$ to the partitions of $n$ given by $\s(\omega w) := \s(w)$ for $\omega 
\in \Omega$.

The partitions of $n$ are naturally partially ordered by 
the dominance order: if $\l$ and $\mu$ are two partitions of $n$, we say that
$\l \domeq \mu$ ($\l$ dominates $\mu$) if for all $k$ we have $$
\sum_{i = 1}^k \l_i \geq \sum_{i = 1}^k \mu_i
.$$ 

\proclaim{Theorem 3.2.4 (Shi)}
Let $(W, S)$ be the extended Coxeter group $\Omega\ltimes W'$ of Definition
3.1.2 corresponding to the affine Hecke algebra of $GL_n$, and let 
$y, w \in W$.  Then $y \leq_{LR} w$ in the
sense of Kazhdan--Lusztig if and only if $\s(y) \domeq \s(w)$.  In
particular, $y$ and $w$ are in the same two-sided cell (\idest $y
\leq_{LR} w \leq_{LR} y$) if and only if $\s(y) = \s(w)$.
\endproclaim

\demo{Proof}
The case where $y, w \in W'$ is dealt with in \cite{{\bf 15}, \S2.9}, and the 
general case is immediate from Lemma 3.2.2.
\qed\enddemo

The relevance of Theorem 3.2.4 for our purposes is that the partitions of
$n$ partially ordered by dominance will form the poset $\Lambda$ of axiom
(A1) for a tabular algebra.

\subhead 3.3 The asymptotic affine Hecke algebra of $GL_n$ \endsubhead

Throughout \S3.3, we shall consider the algebra $(A, \BB)$, where $A$
is the affine Hecke algebra of $GL_n$ from Definition 3.1.2 and $\BB$ is its
Kazhdan--Lusztig basis.
We shall denote the Coxeter group and the extended Coxeter 
group corresponding to $(A, \BB)$ by $(W', S)$ and $(W, S)$ respectively.

The Hecke algebra has an asymptotic algebra that arises
from a construction analogous to that of Definition 2.1.1.  
We recall the definition from \cite{{\bf 17}, \S1.5}.

\definition{Definition 3.3.1}
Let $(A, \BB)$ be the affine Hecke algebra of $GL_n$ equipped with the 
Kazhdan--Lusztig basis.  Define the structure constants $g_{x, y, z}$ by the
formula $$
C_x C_y = \sum_{z \in W} g_{x, y, z} C_z
$$ and define the corresponding integers $\g_{x, y, z}$ as in Definition 1.1.2.
The {\it asymptotic affine Hecke algebra of $GL_n$}, $J$, is the free 
$\zed$-algebra
with basis $\{t_w : w \in W\}$ and structure constants $$
t_x t_y = \sum_{z \in W} \g_{x, y, z} t_z
.$$  If $\bc$ is a two-sided cell in $W$, we define the ring $J_\bc$ to be
the $\zed$-module with free basis $\{t_w : w \in \bc\}$.  
\enddefinition

\proclaim{Lemma 3.3.2}
The algebra $J$ of Definition 3.3.1 decomposes as a direct sum of two-sided 
ideals $J \cong \bigoplus_\bc J_\bc$, where the sum ranges over all two-sided 
cells of $W$.
\endproclaim

\demo{Proof}
This is a consequence of the theory of cells in affine Weyl groups developed
in \cite{{\bf 10}}; see \cite{{\bf 17}, \S1.5} for further discussion of this
property.
\qed\enddemo

\definition{Definition 3.3.3}
A {\it based ring} is a pair $(R, B)$, where $R$ is a unital $\zed$-algebra
with free $\zed$-basis $B$ and nonnegative structure constants.
An isomorphism of based rings is an isomorphism abstract 
$\zed$-algebras compatible with the distinguished bases.
\enddefinition

The asymptotic Hecke algebra $J$ is an example of a based ring; the positivity
property follows from results in \cite{{\bf 9}, \S3} (see \cite{{\bf 17}, \S1.3} for 
the case of extended Coxeter groups).

The structure of the asymptotic algebra of Definition 3.3.1 as a based
ring was established by Xi in \cite{{\bf 17}}, thus verifying a conjecture of 
Lusztig \cite{{\bf 12}, Conjecture 10.5} in a the case of type $A$.  Lemma 3.3.2 
quickly reduces this problem to that of understanding the algebras $J_\bc$.

\proclaim{Theorem 3.3.4 (Xi)}
Let $\bc$ be the two-sided cell associated by the map $\s$ of Theorem 3.2.4
to the partition $\l$ of $n$.  Then $J_\bc$ is isomorphic as a based ring
to the full matrix ring $M_{n_\mu}(\Gamma(\l))$ over a certain (commutative)
based ring $\Gamma(\l)$.  Here, $n_\mu = n!/(\mu_1! \cdots \mu_r!)$, where 
$\mu$ is the dual partition of $\lambda$.
\endproclaim

\demo{Proof}
This follows from \cite{{\bf 17}, Theorem 2.3.2}.
\qed\enddemo

\proclaim{Proposition 3.3.5}
Let $\l$ be a partition of $n$, let $\bc$ be the two-sided cell corresponding 
to $\l$ as in Theorem 3.2.4 and let $L$ be a left cell contained in $\bc$.
Let $J_{L \cap L^{-1}}$ be the $\zed$-module with basis 
$\{t_w : w \in L \cap L^{-1}\}$.  Then $J_{L \cap L^{-1}}$ is a subring of $J$
isomorphic (as a based ring) to the based ring $\Gamma(\l)$ of Theorem 3.3.4.
Furthermore, these based rings are normalized table algebras over $\zed$,
and the table algebra anti-automorphism of $J_{L \cap L^{-1}}$ is induced
by the map $t_w \mapsto t_{w^{-1}}$.
\endproclaim

\demo{Proof}
It follows from \cite{{\bf 10}, Theorem
1.10} and Lemma 3.2.2 that $L$ contains a unique distinguished involution, $d$,
in the sense of Definition 3.1.6.  It follows easily from Definition
3.2.1 that $L^{-1}$ is a right cell; it also contains $d$ since $d^2 = 1$.

It can be shown using standard properties of asymptotic Hecke algebras that
$J_{L \cap L^{-1}}$ is a subring of $J$ (see \cite{{\bf 17}, \S1.5}),
and this subring is isomorphic as a based ring to $\Gamma(\l)$ by 
\cite{{\bf 17}, Theorem 2.3.2}.  We now show that $J_{L \cap
L^{-1}}$ is a normalized table algebra.

The identity element is $t_d$, which follows from the properties of the
integers $\g$ developed in \cite{{\bf 10}}; see \cite{{\bf 17}, \S1.5} for remarks on
the extended Coxeter group case.

Axiom (T1) follows from the positivity properties of \cite{{\bf 9}, \S3}; the
structure constants are integers by definition.

The map $\bar{\ } : t_w \mapsto t_{w^{-1}}$ clearly leaves the set 
$L \cap L^{-1}$ stable and permutes the basis of $\Gamma(\l)$.  
Since $\afn(z) = \afn(z^{-1})$ by \cite{{\bf 9}, Proposition 2.2}, we can
apply the map $*$ of Lemma 3.1.5 to a product $C_x C_y$ to show that
that $g_{x, y, z} = g_{y^{-1}, x^{-1}, z^{-1}}$ and $\g_{x, y, z} =
\g_{y^{-1}, x^{-1}, z^{-1}}$, which proves that $\bar{\ }$ is an 
anti-automorphism satisfying axiom (T2).
(Compare with \cite{{\bf 10}, 1.1 (f)}, which proves the above identity for $\g$
in the case of Coxeter groups.)

Another standard property of $\g$ that follows from the results of \cite{{\bf 10}}
is $\g_{x, y, z} = \g_{y, z^{-1}, x^{-1}}$ (see \cite{{\bf 17}, \S1.3 (b)}).  
Combining this with the 
property in the previous paragraph, we find that $\g_{x, y, z} = 
\g_{z, y^{-1}, x}$.  This is precisely the property needed to prove axiom (T3).
\qed\enddemo

\subhead 3.4 Tabular structure of the affine Hecke algebra of $GL_n$ 
\endsubhead

We continue to concentrate on $(A, \BB)$, the affine Hecke algebra of $GL_n$
(equipped with the Kazhdan--Lusztig basis),
and on its extended Coxeter group $(W, S)$.  
To prove axiom (A3), we need to take a closer look at the structure constants
arising from the basis $\BB$.

Definitions 3.4.1 and 3.4.2 appear in \cite{{\bf 13}, \S1}.

\definition{Definition 3.4.1}
Let $(A, \BB)$ be as above, and let $\bc$ be a two-sided Kazhdan--Lusztig 
cell.  We define
$A_\bc$ to be the free $\A$-module with $\bc$ as a basis.  This has the
structure of an $A$--$A$ bimodule: the left module structure is
given by the formula $$
b . b' = \sum_{b'' \in \bc} g_{b, b', b''} b''
,$$ where $b \in \BB$ and $b' \in \bc$, and the right module structure is
given by the same formula but with $b \in \bc$ and $b' \in \BB$.  The two
module structures commute by associativity of $A$ and the partial order on 
the cells.
\enddefinition

We now introduce a second indeterminate, $v'$.  We denote by $A'$ the
$\zed[v', v^{\prime -1}]$-algebra obtained from $A$ by substituting $v'$ 
for $v$.  We write $g_{b, b', b''}(v)$ for $g_{b, b', b''}$ to emphasize that
$g_{b, b', b''} \in \A$.

\definition{Definition 3.4.2}

Maintain the above notation.
Let $G_\bc$ be the free $\zed[v, v^{-1}, v', v^{\prime -1}]$-module with basis
$\bc$.  We endow $G_\bc$ with the structure of a left $A$-module using the 
formula in Definition 3.4.1.  We also endow $G_\bc$ with the structure
of a right $A'$-module, also using the formula in Definition 3.4.1 but
substituting $g_{b, b', b''}(v')$ for $g_{b, b', b''}(v)$.
\enddefinition

\proclaim{Proposition 3.4.3}
The module $G_\bc$ of Definition 3.4.2 is an $A$--$A'$ bimodule, \idest the
two module structures commute.
\endproclaim

\demo{Proof}
It is enough to check that expressions of the form $b . b' . b''$ denote well 
defined elements of $G_\bc$ whenever $b$ is a basis element of $A$, $b''$ is
a basis element of $A'$ and $b'$ is a basis element of $G_\bc$.  Let $w, 
w', w''$ denote the elements of $W$ indexing the basis elements $b, b', b''$
respectively.

Remark 3.1.3 reduces the problem to the consideration of the case where
$w$, $w'$, $w''$ all lie in the Coxeter group $W'$, because computation of the
product $C_\omega C_u$ does not involve $v$ when $\omega \in \Omega$.
The proof is completed by \cite{{\bf 13}, Theorem 2.2}.
\qed\enddemo

\proclaim{Proposition 3.4.4}
Maintain the notation of Definition 3.4.1.
For any $b_1$, $b_2$, $b_3$, $\be' \in \BB$ such that $\be' \in \bc$ and 
$b_2 \in \bc$, we have $$
\sum_{\be \in \bc} g_{b_1, b_2, \be}(v) \g_{\be, b_3, \be'}
= 
\sum_{\be \in \bc} g_{b_1, \be, \be'}(v) \g_{b_2, b_3, \be}
.$$
\endproclaim

\demo{Proof}
This is formally the same as \cite{{\bf 13}, Proposition 1.9 (a)}.  The proof is
valid in this context due to Proposition 3.4.3 and the fact (proved in
\cite{{\bf 9}, Theorem 5.4}) that the $\afn$-function is constant on elements 
of $\bc$.
\qed\enddemo

We now have all the ingredients to prove the main result.

\proclaim{Theorem 3.4.5}
The Kazhdan--Lusztig basis is a tabular basis for the affine Hecke algebra
of $GL_n$ ($n \geq 3$), and the distinguished involutions of the tabular
basis agree with the distinguished involutions in the sense of Lusztig.
\endproclaim

\demo{Proof}
Let $\Lambda$ be the set of partitions of $n$, partially ordered by dominance
so that $\l \leq \mu$ means $\l \domeq \mu$.
For each $\l \in \Lambda$, let $M(\l)$ be the set of numbers from $1$ up to
$n_\mu$, where $n_\mu$ is as in Theorem 3.3.4, and let $(\Gamma(\l), B(\l))$
be the normalized table algebra of Proposition 3.3.5.  The map $C$ of the
table datum takes a triple $(m, b, m')$ from $M(\l) \times B(\l) \times M(\l)$
and associates to it the Kazhdan--Lusztig basis element $C_w$ (where $w \in W$)
for which $t_w$ corresponds to the element $e_{m m'} \otimes b$ under the
isomorphism of Theorem 3.3.4.  The map $*$ sends $C_w$ to $C_{w^{-1}}$.

We now check the axioms of Definition 1.1.3.

Axiom (A1) follows easily once we notice that the identity element of the
Hecke algebra is a Kazhdan--Lusztig basis element.

We saw in Lemma 3.1.5 that $*$ is an anti-automorphism.
We now need to check that it is compatible with the map $C$ of the table
datum.
Fix $\l \in \Lambda$, let $\bc$ be the corresponding two-sided cell 
(as in Theorem 3.2.4) and let $\bar{\ }$ denote the table algebra 
anti-automorphism of $\Gamma(\l)$ as in Proposition 3.3.5.
We saw in the proof of that proposition that the map $t_w \mapsto t_{w^{-1}}$
extends to an algebra anti-automorphism of $J$.
It follows from Theorem 3.3.4, Proposition 3.3.5 and \cite{{\bf 17}, Lemma 2.3.1} 
that this anti-automorphism acts on $J_\bc$ by sending 
$e_{m m'} \otimes b$ to $e_{m' m} \otimes \bar{b}$,
where $\bar{\ }$ is the table algebra anti-automorphism of $\Gamma(\l)$.
Axiom (A2) now follows.

The claims about the partial order in axiom (A3) follow from Lemma 3.2.3 and
Theorem 3.2.4.  
For the other claims, we note that the left $A$-module $A_\bc$ of Definition 
3.4.1 is a right $J_\bc$-module affording the regular representation after
the two bases are identified in the obvious way.  Proposition 3.4.4 says that 
these two actions commute with each other.  Now label each basis element of
$A_\bc$ by an element $e_{m m'} \otimes b$ as given by the isomorphism of
Theorem 3.3.4, and the remaining claims are consequences of Proposition 3.4.4.

Axiom (A4) follows easily from Theorem 3.3.4, which allows the elements
\newline $\g_{K, K', K''}$ to be computed.  The axiom says that $$
(e_{m_1 m_2} \otimes b_1)(e_{m_3 m_4} \otimes b_2)
$$ contains the basis element $e_{m_5 m_6} \otimes b_3$ with nonzero 
coefficient if and only if both factors come from the same summand $J_\bc$,
$m_5 = m_1$, $m_2 = m_3$, $m_6 = m_4$ and $b_3$
occurs with nonzero coefficient in $b_1 b_2$. The assertion regarding
$\g_{K, K', K''} = 1$ is also clear.

For axiom (A5), we define an $\A$-linear function 
$\t$ on the Hecke algebra $\H$ of $(W, S)$ by the condition $$
\t \left( \sum a_w T_w \right) := a_1
.$$  We claim $\t(x y) = \t(yx)$ for all $x, y \in \H$; it is enough to check
that this holds when $x = T_s$ or $x = T_\omega$ for $\omega \in \Omega$,
and this follows easily from the relations in \S3.1.  (The corresponding
result for Coxeter groups is \cite{{\bf 9}, 1.4.1}.)  Since $*$ sends $T_w$ to
$T_{w^{-1}}$, it is clear that $\t(x) = \t(x^*)$.  

To complete the proof, let us denote by $(W', S)$ the Coxeter group 
corresponding to $(W, S)$.  If $z \in W'$, it follows from Definition 
3.1.6 and the formula for $C_z$ in Definition 3.1.1 that $$
\t(v^{\afn(C_z)} C_z) = 
\cases a \mod v^{-1} \A^- & \text{ if } z \in {\Cal D}, \cr
0 \mod v^{-1} \A^- & \text{ otherwise,} \cr
\endcases
$$ where $a \in \zed$.  By \cite{{\bf 10}, Proposition 1.4 (a)}, we see that
$P_{e, z}$ is monic when $z \in {\Cal D}$, and this proves that $a = 1$.

Now if $u \in \Omega$, Remark 3.1.3 shows that $\t(C_{uz})$ is zero unless
$u$ is the identity, and we have already dealt with the latter case.  We 
conclude that $$
\t(v^{\afn(C_w)} C_w) = 
\cases 1 \mod v^{-1} \A^- & \text{ if } z \in {\Cal D}, \cr
0 \mod v^{-1} \A^- & \text{ otherwise.} \cr
\endcases
.$$

It remains to show that the distinguished involutions in the sense of 
Definition 3.1.6 coincide with the distinguished involutions of axiom (A5).
As mentioned in \cite{{\bf 17}, \S1.5}, the ring $J$ has an identity, namely
$\sum_{d \in \D} t_d.$  From the construction in \cite{{\bf 4}, Lemma 3.2.2}, we 
now see that $$
\{C_d : d \in {\Cal D}\} = \{C_{S, S}^1 : S \in M(\l), \l \in \Lambda\}
,$$ where the $C$ on the right is as in the definition of table datum.  
This means that the two senses of the term ``distinguished involution'' agree
in this case, and the proof is complete.
\qed\enddemo

\head 4. Applications \endhead

In \S4, we compare the results of this paper with some related results 
in the literature.

\subhead 4.1 Description of Kazhdan--Lusztig cells using the table datum 
\endsubhead

It is not hard to show that the cell modules for the affine Hecke 
algebra of $GL_n$ as a tabular algebra are compatible with the
cell representations in the sense of Kazhdan--Lusztig; we do this now for
the sake of easy reference.

\proclaim{Proposition 4.1.1}
Let $A$ be the affine Hecke algebra of $GL_n$ ($n \geq 3$), let $\BB$ be its 
Kazhdan--Lusztig basis, and let $(W, S)$ be the corresponding Coxeter
system.  Let $(\Lambda, \Gamma, B, M, C, *)$ be the table
datum for $(A, \BB)$ given in the proof of Theorem 3.4.5.  Let 
$x, y \in W$, and write $C_{T, U}^b = C_x$ and $C_{V, W}^{b'} = C_y$, 
where $C_x, C_y \in \BB$, $T, U \in M(\l)$ and $V, W \in M(\l')$ for 
some $\l, \l' \in \Lambda$.
Then:
\item{\rm (i)}{$x \sim_{LR} y \Leftrightarrow \l = \l'$;}
\item{\rm (ii)}{$x \sim_L y \Leftrightarrow \l = \l'$ and $U = W$;}
\item{\rm (iii)}{$x \sim_R y \Leftrightarrow \l = \l'$ and $T = V$.}
\endproclaim

\demo{Proof}
We first prove (ii).  Two elements $C_{S, T}^b$ and $C_{U, V}^{b'}$ 
are in the same left cell of $\BB$ (in
the sense of tabular algebras given in Lemma 1.2.5) if and only if $\l = \l'$
and $U = W$.  The preorder generating this equivalence relation is given in
Definition 1.2.3, and the latter is clearly compatible with the 
characterization of the preorder $\leq_L$ given in Lemma 3.2.3.  The proof
of (ii) now follows from the definition of left cells in the Kazhdan--Lusztig
sense.

Part (iii) follows from (ii) by applying the tabular anti-automorphism $*$.
This sends $C_{S, T}^b$ to $C_{T, S}^{\bar{b}}$ by definition, and sends
$C_w$ to $C_{w^{-1}}$ by Lemma 3.1.5; recall from Definition 3.2.1 that
inversion sends left cells to right cells (in the sense of Kazhdan--Lusztig).

The proof of (i) is similar to the proof of (ii), replacing Lemma 1.2.5 by
\cite{{\bf 4}, Proposition 3.1.3}, Definition 1.2.3 by \cite{{\bf 4}, Definition 3.1.1}
and using the characterization of $\leq_{LR}$ in Lemma 3.2.3.
\qed\enddemo

\remark{Remark 4.1.2}
Although there are many possible table data for the pair $(A, \BB)$, 
Proposition 4.1.1 can be shown to be independent of the choice of table
datum using ideas from \cite{{\bf 6}}.
\endremark

One application of Proposition 4.1.1 is that it shows that, in the situation
under consideration, the cell modules of \S1.2 afford the left cell 
representations defined in \cite{{\bf 8}}.

\subhead 4.2 Standard modules for affine Hecke algebras \endsubhead

Combining Theorem 2.2.3 with Theorem 3.4.5 gives a construction of a set
of ``standard modules'' for the affine Hecke algebra of $GL_n$.  (Notice
that, by Definition 1.3.4 and Theorem 2.2.3 (ii), these modules can be 
defined in terms of the
representation theory of the table algebras $\Gamma(\l)$ without reference to
the asymptotic algebra.)  The next result justifies the use of the term
``standard modules'' introduced in \S1.3.

\proclaim{Proposition 4.2.1}
Let $A$ be the affine Hecke algebra of $GL_n$ ($n \geq 3$), let $\BB$ be its 
Kazhdan--Lusztig basis, and let $(W, S)$ be the corresponding Coxeter
system.  The (isomorphism classes of) standard modules associated to the 
pair $(A, \BB)$ in the sense of \S1.3 agree with the standard modules
${\Cal K}_{u, s, \rho}$ in the sense of Lusztig \cite{{\bf 12}}.
\endproclaim

\demo{Note}
Recall from Theorem 2.2.3 (iv) that the isomorphism classes of standard
modules for a tabular algebra depend only on the basis.
\enddemo

\demo{Proof}
The modules ${\Cal K}_{u, s, \rho}$ are defined in \cite{{\bf 12}}, and are called
``standard modules'' in the introduction to \cite{{\bf 12}}.  By 
\cite{{\bf 12}, Theorem 4.2}, each module ${\Cal K}_{u, s, \rho}$ is isomorphic
to precisely one module $^\Phi E$, where $E$ is a simple $J$-module made
into an $\H$-module by applying a certain homomorphism $\Phi : \H \ra J$.
(The set-up in \cite{{\bf 12}} is that the base ring is a field, $k$, and the 
parameter $v$ acts by scalar multiplication by $r \in k^*$, as in our \S2.2.)
The map $\Phi$ is defined to be $\bigoplus_{\l \in \Lambda} \Phi_\l$, where
$\Phi_\l$ is as in Lemma 2.1.4 and $\Lambda$ is as in the proof of Theorem
3.4.5; this is a restatement of the formula in
\cite{{\bf 12}, \S1.4} using the identifications of Theorem 3.4.5 and Proposition
4.1.1.  Since $A_\l^\infty$ (\idest $J_\bc$ for some two-sided cell $\bc$)
is a direct summand of $J$ and $E$ is simple, $^\Phi E$ is isomorphic to
$^{\Phi_\l} E$ for some $\l \in \Lambda$.  The conclusion follows from
Theorem 2.2.3 (iii).
\qed\enddemo

\remark{Remark 4.2.2}
An interesting problem is the determination of the decomposition matrix of
the standard modules into simple modules.  If the scalar $r$ (as in the 
proof above) is not a root of unity or is equal to $1$, \cite{{\bf 11}, Theorem 3.4}
shows that the standard modules are simple, but at a root of unity the
situation is much more complicated.
\endremark

\subhead 4.3 Relationship with two-step nilpotent representations \endsubhead

We now outline how we can use the language of tabular algebras to set up some 
of the main results in Graham and Lehrer's work on two-step nilpotent 
representations of the Hecke algebra of $GL_n$ \cite{{\bf 3}}.

As usual, we let $A$ be the affine Hecke algebra of $GL_n$ ($n \geq 3$), 
let $\BB$ be its Kazhdan--Lusztig basis, and let $(W, S)$ be the 
corresponding Coxeter system.
It follows (for example from \cite{{\bf 5}, Theorem 3.4}) that if we consider the 
set $$\Lambda_c := \{\l \in \Lambda: \l_1 \geq 3\}$$ of all partitions of $n$
whose first part exceeds $3$, then the set of basis elements 
$\bigcup_{\l \in \Lambda_c} \bc_\l$ spans an ideal of $A$.  Quotienting $A$
by this ideal yields an algebra denoted by $\widetilde{TL}_n^a(q)$ in 
\cite{{\bf 3}, (1.7)}.  It is an easy consequence of Theorem 3.4.5 that this
algebra is tabular with tabular basis given by the nonzero images of the
elements $C_w$.  This basis is closely related to, but not always the same 
as, the basis used in \cite{{\bf 3}}.  Either basis may be handled using only
combinatorics, which is not the case for the full Kazhdan--Lusztig basis $\BB$.

The ``cell modules'' of \cite{{\bf 3}} are simply
the standard modules of $\widetilde{TL}_n^a(q)$ with the projection of the
Kazhdan--Lusztig basis as the tabular basis.
(The construction of these modules is rather similar
to our Definition 1.3.4.)  One of the main achievements of \cite{{\bf 3}} is the
determination of the decomposition matrix of the standard modules into
simples.  It turns out (see, for example, \cite{{\bf 3}, Theorem 5.5 (ii)})
that each standard module has a unique simple quotient, and every simple
module is isomorphic to one arising in this way.  This behaviour is familiar
from the theory of cellular algebras in \cite{{\bf 2}}.

\subhead 4.4 Concluding remarks \endsubhead

As well as proving Lusztig's conjecture for the affine Hecke algebra of $GL_n$,
Xi also considers the case where $\Omega$ is the group $\zed$, also acting
by rotations of the Coxeter graph.  The conclusion (given in \cite{{\bf 17}, \S8.2})
is that the conjecture holds, and we have an analogue of Theorem 3.3.4.
A version of Theorem 3.4.5 will also hold, with some light modifications to
the proof.

Xi's main results in \cite{{\bf 17}} may be paraphrased in the language of tabular 
algebras by saying that they verify that the table data of certain extended
affine Hecke algebras are as predicted by Lusztig; this is a highly
nontrivial result.  An independent proof of Theorem 3.4.5 would achieve some of
the steps required in the proofs of the main results of \cite{{\bf 17}, \S8.4}.

As regards standard modules, it would be very interesting to know whether
the behaviour described in \S4.3 is typical; that is, can one always construct
a set of simple modules from the heads of the standard modules?  One
problem to be solved here is a replacement for the bilinear form used in
\cite{{\bf 2}} or \cite{{\bf 3}} to construct the simple modules from the standard
modules, especially in the case where the table algebras involved are
noncommutative.

\leftheadtext{}
\rightheadtext{}

\end
\Refs\refstyle{A}\widestnumber\key{{\bf 10}}
\leftheadtext{References}
\rightheadtext{References}

\ref\key{{\bf 1}}
\by Z. Arad and H.I. Blau
\paper On Table Algebras and Applications to Finite Group Theory
\jour J. Algebra
\vol 138 \yr 1991 \pages 137--185
\endref

\ref\key{{\bf 2}}
\by J.J. Graham and G.I. Lehrer
\paper Cellular algebras
\jour Invent. Math.
\vol 123
\yr 1996
\pages 1--34
\endref

\ref\key{{\bf 3}}
\by J.J. Graham and G.I. Lehrer
\paper The two-step nilpotent representations of the extended affine
Hecke algebra of type $A$
\jour Compositio Math.
\toappear
\endref

\ref\key{{\bf 4}}
\by R.M. Green
\paper Tabular algebras and their asymptotic versions
\jour J. Algebra
\vol 252 \yr 2002 \pages 27--64
\endref

\ref\key{{\bf 5}}
\by R.M. Green
\paper On $321$-avoiding permutations in affine Weyl groups
\jour J. Algebraic Combin.
\vol 15 \yr 2002 \pages 241--252
\endref

\ref\key{{\bf 6}}
\by R.M. Green
\paper Categories arising from tabular algebras
\jour Glasgow Math. J.
\miscnote to appear; \newline {\tt math.QA/0207097}
\endref

\ref\key{{\bf 7}}
\by R.M. Green and J. Losonczy
\paper Fully commutative Kazhdan--Lusztig cells
\jour Ann. Inst. Fourier (Grenoble)
\vol 51 \yr 2001 \pages 1025--1045
\endref

\ref\key{{\bf 8}}
\by D. Kazhdan and G. Lusztig
\paper Representations of Coxeter groups and Hecke algebras
\jour Invent. Math. 
\vol 53 \yr 1979 \pages 165--184
\endref

\ref\key{{\bf 9}}
\by G. Lusztig
\paper Cells in affine Weyl groups
\inbook Algebraic groups and related topics
\publ Adv. Studies Pure Math 6, North-Holland and Kinokuniya
\publaddr Tokyo and Amsterdam
\yr 1985
\pages 255--287
\endref

\ref\key{{\bf 10}}
\by G. Lusztig
\paper Cells in affine Weyl groups, II
\jour J. Alg.
\vol 109
\yr 1987
\pages 536--548
\endref

\ref\key{{\bf 11}}
\by G. Lusztig
\paper Cells in affine Weyl groups, III
\jour J. Fac. Sci. Tokyo U. (IA)
\vol 34
\yr 1987
\pages 223--243
\endref

\ref\key{{\bf 12}}
\by G. Lusztig
\paper Cells in affine Weyl groups, IV
\jour J. Fac. Sci, Tokyo U. (IA)
\vol 36
\yr 1989
\pages 297--328
\endref

\ref\key{{\bf 13}}
\by G. Lusztig
\paper Quantum groups at $v = \infty$
\jour Prog. Math.
\vol 131 \yr 1995 \pages 199--221 
\endref

\ref\key{{\bf 14}}
\by J.Y. Shi
\paper The Kazhdan--Lusztig cells in certain affine Weyl groups
\jour Lecture Notes in Mathematics
\vol 1179
\yr 1986
\publ Spinger
\publaddr Berlin
\endref

\ref\key{{\bf 15}}
\by J.Y. Shi
\paper The partial order on two-sided cells of certain affine Weyl
groups
\jour J. Alg.
\vol 176
\yr 1996
\pages 607--621
\endref

\ref\key{{\bf 16}}
\by V.S. Sunder
\paper $\text{II}_1$ factors, their bimodules and hypergroups
\jour Trans. Amer. Math. Soc.
\vol 330 \yr 1992 \pages 227--256
\endref

\ref\key{{\bf 17}}
\by N. Xi
\paper The based ring of two-sided cells of affine Weyl groups of type 
$\tilde A_{n-1}$
\jour Mem. Amer. Math. Soc.
\vol 157 \yr 2002 \pages no. 749
\endref

\endRefs
